\documentclass[12pt]{amsart}

\usepackage{amssymb,amscd,pb-diagram,mathrsfs}
\usepackage{fullpage}

\def\C{{ \mathbb C }}

\def\Z{{   \mathbb Z }}
\def\Q{{\mathbb Q}}

\def\F{{\mathbb F}}
\def\1{{\bf 1}}

\def\NN{{ \mathscr N}}

\def\ord{{\text{ord}}}

\def\P{{\mathbb P}}

\def\l{{\ell}}

\def\X{ {\mathcal{X}}}
\def\Y{ {\mathcal{Y}}}

\def\c2{\chi_2}

\def\hfl#1#2{\smash{\mathop{\hbox to
12mm{\rightarrowfill}}\limits_{\scriptstyle #1}^{\scriptstyle #2}}}

\def\C{{ \mathbb C }}

\def\Z{{   \mathbb Z }}
\def\Q{{\mathbb Q}}

\def\F{{\mathbb F}}

\def\P{{\mathbb P}}

\def\l{\ell}

\def\deg{{ \mathop{{\rm {deg}}} }}

\def\iL{ {\mathcal L}_n^{(\alpha)}(x) }
\def\IL{ {\mathcal L}_n^{\langle r \rangle}(x)}

\def\LL{L_n^{\langle r \rangle}(x)}
\def\L#1{L_n^{\langle #1\rangle}(x)}

\newcount\timehh\newcount\timemm
\timehh=\time 
\divide\timehh by 60 \timemm=\time
\newtheorem{theorem}{Theorem}[section]
\newtheorem{lemma}[theorem]{Lemma}
\newtheorem{corollary}[theorem]{Corollary}

\newtheorem{conjecture}[theorem]{Conjecture}
\newtheorem{definition}[theorem]{Definition}

\newtheorem{question}[theorem]{Question}

\newenvironment{Proof}{\removelastskip\par\medskip
\noindent{\em Proof.} \rm}{\penalty-20\null\hfill$\square$\par\medbreak}

\title[Algebraic Properties of Laguerre Polynomials]
{Algebraic Properties of A Family of Generalized Laguerre Polynomials}
\author{
          Farshid Hajir\\
   {
    \protect \protect\sc\today\ -- 
    \ifnum\timehh<10 0\fi\number\timehh\,:\,\ifnum\timemm<10 0\fi\number\timemm
    \protect \, \, \protect 
  }
}
\thanks{This work was supported by the National Science
Foundation under Grant No. 0226869.}
\begin{document}

\begin{abstract}
We study the algebraic properties of Generalized Laguerre Polynomials
for negative integral values of the parameter.  For integers $r,n\geq
0$, we conjecture that $L_n^{(-1-n-r)}(x) = \sum_{j=0}^n
\binom{n-j+r}{n-j}x^j/j!$ is a $\Q$-irreducible polynomial whose
Galois group contains the alternating group on $n$ letters.  That this
is so for $r=n$ was conjectured in the 50's by Grosswald and proven
recently by Filaseta and Trifonov.  It follows from recent work of
Hajir and Wong that the conjecture is true when $r$ is large with
respect to $n\geq 5$.  Here we verify it in three situations: i) when
$n$ is large with respect to $r$, ii) when $r \leq 8$, and iii) when
$n\leq 4$.  The main tool is the theory of $p$-adic Newton Polygons.
\end{abstract}

\maketitle

\section{Background and Summary of Results}

The Generalized Laguerre Polynomial (GLP) is a one-parameter family
defined by 
$$L_n^{(\alpha)}(x) = (-1)^n\sum_{j=0}^n \binom{n+\alpha}{n-j}
\frac{(-x)^j}{j!}.$$ Here, as usual, the binomial coefficient
$\binom{t}{k}$ is defined to be $t(t-1)\cdots(t-k+1)/k!$ for
non-negative integers $k$; the inclusion of the sign $(-1)^n$
is not standard.  Sometimes it is more convenient to work
with the monic integral polynomial $\iL=n!L_n^{(\alpha)}(x)$.  The
monographs by P\'olya-Szeg\H o \cite{pz}, Szeg\H o \cite{op}, and
Andrews-Askey-Roy \cite{aar} contain a wealth of facts about this and
other families of orthogonal polynomials.
To cite only two,
we have the second order linear (hypergeometric) differential equation
$$
x y'' + (\alpha + 1 -x) y' + n y = 0, \qquad y= L_n^{(\alpha)}(x),
$$
as well as the difference equation
$$ L_n^{(\alpha-1)}(x) - L_n^{(\alpha)}(x)=L_{n-1}^{(\alpha)}(x).
$$ 
A quick glance at the mathematical literature makes it clear that GLP
has been extensively studied primarily because of the very important
roles it plays in various branches of analysis and mathematical
physics.  However, not long after its appearance in the literature early
in the twentieth century, it became evident, in the hands of Schur,
that GLP also enjoys {\em algebraic} properties of great interest.


For instance, in 1931, Schur \cite{sch70} gave a pretty formula for the
discriminant of ${\mathcal L}_n^{(\alpha)}(x)$:
\begin{equation}
\label{disc} \Delta_n^{(\alpha)} =
\prod_{j=2}^n j^j (\alpha+j)^{j-1}.
\end{equation}
In \cite{sch67} and \cite{sch70}, he showed that $L_n^{(0)}(x)$
(classical Laguerre polynomial, first studied by Abel), and
$L_n^{(1)}(x)$ (derivative of classical Laguerre), are irreducible in
$\Q[x]$ for all $n$; he also calculated their Galois groups.

Recently, a number of articles concentrating on the algebraic
properties of GLP have appeared, including Feit \cite{feit}, Coleman
\cite{coleman}, Gow \cite{gow}, Hajir \cite{jnt},
Filaseta-Williams\cite{fw}, Sell \cite{wrobel}.  In all of these papers, the
authors take a sequence $(\alpha_n)_n$ of rational numbers and
consider the irreducibility and Galois group of $L_n^{(\alpha_n)}(x)$
over $\Q$.  The best general such result to date is for
constant sequences $\alpha_n$.  
\par
\vspace{\baselineskip}
\noindent {\bf Theorem.} (Filaseta-Lam/Hajir) {\em Suppose $\alpha$ is a
fixed rational number which is not a negative integer.  Then
for all
but finitely many integers $n\geq 0$, $L_n^{(\alpha)}(x)$ is irreducible over
$\Q$ and has Galois group containing $A_n$.}
\par
\vspace{\baselineskip}
\noindent It should be noted that reducible GLP for rational values of
the parameter $\alpha$ do exist (already infinitely many exist in
degrees 2, 3 or 4, cf. Section 6).  The irreducibility part of the
above theorem is due to Filaseta and Lam \cite{fl}; the supplement on
the Galois group was added in \cite{hajir}.  The proof of both parts
is effective.



At the values of the parameter $\alpha$ excluded by the theorem of
Filaseta and Lam (the negative integers), one finds some of the most
interesting families of GLP, e.g.  the truncated exponential series,
and the Bessel Polynomials (see below).  In this paper, we consider
irreducibility and Galois groups of GLP for exactly these values of
the parameter $\alpha$.  Note that their exclusion from the theorem is
quite necessary; namely, when $\alpha$ is a negative integer,
$L_n^{(\alpha)}(x)$ is reducible for all $n\geq |\alpha|$.  Indeed,
writing $\alpha = -a$ with $n=a+m$ where $a$ is an integer in $[1,n]$
we have
\begin{equation}\label{factor}
{\mathcal L}_n^{(-a)}(x) = x^a\cdot {\mathcal L}_m^{(a)}(x), \qquad
\mathcal{L}_m^{(a)}(0) \neq 0. 
\footnote{Incidentally the repeated roots at
the origin evident in the above factorization (for $2\leq a \leq n$
i.e. $-n\leq \alpha \leq -2$) explain the presence of the factors
$\alpha+j$, $j=2,\ldots, n$, in (\ref{disc}).  Their multiplicites in
the discriminant (i.e. $j-1$) express the tameness of the
corresponding ramified points in the extension 
$\C(\alpha)\hookrightarrow \C(\alpha)[x]/(L_n^{(\alpha)}(x))$ of
function fields.  It would be interesting to obtain a similarly conceptual
explanation of the factors $j^j$ as well.}
\end{equation}

Given the above observation, namely that for small negative integral values
of the parameter $\alpha$, $L_n^{(\alpha)}(x)$ is a simple factor times a
Laguerre polynomial of positive parameter, it is natural to replace the
parameter $\alpha$ by a parameter $r$ via the translation
$$\alpha = -1 - n - r, $$
and to consider instead
\begin{eqnarray}\label{lnr} \nonumber \LL &:=& L_n^{(-1-n-r)}(x)\\
&=& \sum_{j=0}^n \binom{n-j+r}{n-j}\frac{x^j}{j!}.
\end{eqnarray}
It is also useful to note that 
\begin{equation}\label{lnr2}
\IL := n! \LL = \sum_{j=0}^n \binom{n}{j}
(r+1)(r+2)\cdots(r+n-j) x^j,
\end{equation}
is monic and has positive integer coefficients, assuming, as we do throughout
the paper, that $r$ is a non-negative integer.

The parametrization (\ref{lnr}) is a natural one in some respects
(in addition to being a convenient representation of the family of
polynomials we wish to consider).  For instance, differentiation with
respect to $x$ of $L_n^{(\alpha)}(x)$ has the effect of lowering $n$
by $1$ and raising $\alpha$ by $1$, so in the new parametrization,
differentiation leaves $r$ fixed:
$$\partial_x \LL = L_{n-1}^{\langle r\rangle}(x).$$ Indeed, the most
familiar such ``derivative-coherent'' sequence of polynomials, namely
the truncations of the exponential series, is obtained when we set
$r=0$:
$$E_n(x):=\L{0} = \sum_{j=0}^n \frac{x^j}{j!}.$$

Let us review some known algebraic facts about $\LL$ for small $r\geq
0$.  The exponential Taylor polynomails $E_n$ were first studied by
Schur.  He showed that they are irreducible over $\Q$ \cite{sch67},
and have Galois group $A_n$ or $S_n$ (over $\Q$) according to whether
$n$ is divisible by $4$ or not \cite{sch70}.  Coleman \cite{coleman}
gave a different proof of these results.  For the case $r=1$,
irreducibility and the calculation of the Galois group using methods
of Coleman and Schur, respectively, were established in \cite{jnt}.
Moreover, in \cite{jnt}, the values of $n$ for which the splitting
field of $\L{0}$ or $\L{1}$ can be embedded in an
$\tilde{A_n}$-extension were determined using formulae of Feit
\cite{feit} and a criterion of Serre \cite{serre}.  All of the above
was carried out for $r=2$ by Sell in \cite{wrobel}.  But perhaps the
best-studied family of GLP is that of Bessel Polynomials (BP) $z_n(x)$
which are, simply the monic GLP with $r=n$.  Namely we have
$$z_n(x):= \sum_{j=0}^n \frac{(2n-j)!}{j!(n-j)!} x^j = {\mathcal
L}_n^{\langle n\rangle}(x).$$ Grosswald pointed out that the BPs play a
distinguished role among GLPs due to certain ``symmetries'' which in our
notation amounts to their invariance under exchange of $r$ and $n$.   They are
arithmetically interesting as well (for example the prime 2 does not
ramify in the algebra $\Q[x]/(z_n(x))$ despite the presence of many
powers of 2 in the discriminant of $z_n$, cf. (\ref{disc})).  Their
irreducibility was conjectured by Grosswald \cite{grosswald}, who also
showed that their Galois group is always the full symmetric group
(assuming his conjecture).  The irreducibility of all BPs was proved,
first for all but finitely many $n$ by Filaseta \cite{facta}, and
later for all $n$ by Filaseta and Trifonov \cite{ft}.

As an extension of Grosswald's conjecture, we have
\begin{conjecture}\label{conj}
For integers $r,n\geq 0$, $\LL$ is irreducible over
$\Q$.\end{conjecture}
\begin{conjecture}\label{galconj}
For integers $r,n\geq0$, if $\LL$ is irreducible over $\Q$, then
its Galois group over $\Q$ contains the alternating group $A_n$.
\footnote{
Note that once we know the Galois group of a degree $n$
polynomial $f$ contains $A_n$, then it is either $A_n$ or $S_n$
according to whether the discriminant of $f$ is a square or not; the
latter is easily determined for our polynomials using Schur's formula
(\ref{disc}).
  }
\end{conjecture}

There is already a fair bit of evidence for this pair of conjectures.
As described above, they are true for all $n$ if $r=0,1,2$ or $r=n$.
In Sell \cite{wrobel}, it was shown that $\LL$ is irreducible over
$\Q$ if $\gcd(n,r!)=1$; that is already enough to show that for each
fixed $r$, Conjecture \ref{conj} is true for a positive proportion of
integers $n\geq 0$ (this proportion goes to zero quickly with $r$
however).


Our first and main result is
\begin{theorem}\label{mt}
For a fixed $r\geq 0$, all but finitely many $\LL$ are irreducible over
$\Q$ and have Galois group (over $\Q$) containing $A_n$.
\end{theorem}

For a more precise (effective) statement, see Theorems \ref{it} and \ref{gt}.
The irreducibility part of Theorem \ref{mt} is a companion of sorts
for the Filaseta-Lam Theorem.  
As an illustration of the effectivity of our approach, and to
gather more evidence for Conjectures \ref{conj} and \ref{galconj},
we prove
the following theorem.
\begin{theorem}\label{mtg}
If  $0\leq r \leq 8$, then for all $n$, $\LL$ is irreducible and has Galois
group containing $A_n$ over $\Q$.
\end{theorem}

Investigating the irreducibility of $\LL$ for a fixed $n$ and all
large $r$ has a different flavor; the methods we use here give us only
a weak result (see Corollary \ref{cokl}).  In a joint work with
Wong \cite{HW}, using algebro-geometric and group-theoretic techniques,
we prove that for each fixed $n\geq 5$, over
a fixed number field $K$, all but finitely many $L_n^{(\alpha)}(x)$
are irreducible and have Galois group containing $A_n$.  In particular, for
$n\geq 5$, Conjectures \ref{conj} and \ref{galconj} hold for all 
$r$ large enough with respect to $n$.

Here, we complement the above result of \cite{HW} by showing that
Conjectures \ref{conj} and \ref{galconj} hold for all $r\geq 0$ if
$n\leq 4$ (Theorem \ref{n3}).  As for the possibility of verifying
further cases of these conjectures, the methods used by Filaseta and
Trifonov \cite{ft} in proving the irreducibility of $\LL$ for $r=n$
should hopefully yield results in the middle range where $r \approx
n$.

The basic strategy we use for proving irreducibility of $\LL$ was
developed by Sell \cite{wrobel} for the case $r=2$ as an extension
of the proof for $r=1$ given in \cite{jnt}, which was itself an
adaptation of Coleman's proof \cite{coleman} for the case $r=0$.  Here
is a sketch of it.  We fix $r\geq 0$ and suppose $g$ is a proper
divisor, in $\Q[x]$ of $\LL$.  In Step 1, using a criterion of Coleman
\cite{coleman} formalized by Sell \cite{wrobel}, we show that
$\deg(g)$ is divisible by $n_0$, the largest divisor of $n$ which is
co-prime to $\binom{n+r}{r}$.  Then $\deg(g)/n_0$ is at most $r!$ so
is bounded since $r$ is fixed.  In Step 2, thanks to a criterion of
Filaseta \cite{fschur}, we eliminate this bounded number of
possibilities for $\deg(g)/n_0$, giving the desired contradiction.
For Filaseta's criterion to apply, we require the existence of certain
auxiliary primes and this is where we have to assume that $n$ is large
with respect to $r$ so as to apply results from analytic number theory
on the existence of primes in short intervals; these are gathered
together in section 3.

We should point out that the Coleman and Filaseta criteria
are both based on the theory of $p$-adic Newton polygons (which we
review in the next section).  Indeed, the key idea of Step 1
is the simple observation that if $p$ is a prime divisor of
$n$ which does not divide the constant coefficient of $\LL$, then the
$p$-adic Newton polygons of $\LL$ and $E_n$ coincide.  

For the computation of the Galois group, we 
use the criterion described in \cite{hajir}, which was already
implicit in Coleman \cite{coleman} and is also based on Newton
Polygons.

Finally, a bibliographic comment.  In Grosswald's meticulously written
treatise {\em Bessel Polynomials} \cite{grosswald}, he considers not
just the BP $z_n(x)$ but ``Generalized Bessel Polynomials (GBP)''
$z_n(x;a)$ and gives much information about their algebraic and
analytic properties.  The GBP is just a different parametrization
of GLP, as described on p. 36 of \cite{grosswald}.  Therefore, even
though it is not billed as such, Grosswald's book is a rich source of
information about GLP.

{\bf Acknowledgments.} I would like to thank Professors Filaseta and Wong 
for their helpful remarks.

\section{Irreducibility Criteria}

For a prime $p$ and $z \in \Q^*$, we write $\ord_p(z)$ for the
$p$-adic valuation of $z$: $\ord_p(z) = a$ where $z= p^a m/n$ with
integers $m$ and $n$ not divisible by $p$.  It is convenient to put
$\ord_p(0)=\infty$.  We extend the $p$-adic valuation $\ord_p$ to the
algebraic closure $\overline{\Q}_p$ of the $p$-adic completion $\Q_p$
of $\Q$ in the standard way, see Gouvea \cite{gouvea} for example.

For the convenience of the reader, we recall some facts from the
theory of $p$-adic Newton Polygons as well as a useful corollary due
originally to Dumas \cite{dumas} but rediscovered and used in the
context of GLP by Coleman \cite{coleman}.  References include Gouvea
\cite{gouvea}, Amice \cite{amice}, Artin \cite{artin}, and Hensel-Landsberg
\cite{hensel}; the latter is, to the best of my knowledge, where the
general notion of $p$-adic Newton Polygons originated.  An excellent
survey on the applications of Newton Polygons for irreducibility is
Mott \cite{mott}.

The {\em $p$-adic Newton Polygon} (or $p$-Newton polygon) $NP_p(f)$ of
a polynomial $f(x)= \sum_{j=0}^n c_j x^j \in \Q[x]$ is the lower
convex hull of the set of points $$S_p(f) = \{ (j,\ord_p(c_j))| 0 \leq
j \leq n\}.$$ It is the highest polygonal line passing on or below the
points in $S_p(f)$.  The vertices $(x_0,y_0), (x_1, y_1), \cdots
,(x_r,y_r)$, i.e. the points where the slope of the Newton polygon
changes (including the rightmost and leftmost points) are called the
{\em corners} of $NP_p(f)$; their $x$-coordinates ($0=x_0< x_1 <
\cdots < x_r=n$) are the {\em breaks} of $NP_p(f)$.  For the $i$th edge,
joining $(x_{i-1},y_{i-1})$ to $(x_i,y_i)$, we put $$H_i=y_i-y_{i-1}, W_i=
x_{i} - x_{i-1}, m_i= H_i/W_i, d_i = \gcd(H_i,W_i), \qquad i=1,\cdots, r.
$$ We call these quantities, respectively, the {\em height, width,
slope} and {\em multiplicity} of the $i$th edge.  We also put
$w_i=W_i/d_i$, $h_i=H_i/d_i$, so that $w_i$ is the denominator, in
lowest terms, of $m_i=H_i/W_i=h_i/w_i$.  The $i$th edge is made up of
$d_i$ {\em segments} of width $w_i$.
We call the $i$th edge {\em pure} if its multiplicity $d_i$ is $1$.

\begin{theorem}[Main Theorem of Newton Polgyons]
\label{NP} Let $(x_0,y_0),(x_1,y_1),\ldots,(x_r,y_r)$ denote the
successive vertices of $NP_p(f).$ Then there exist polynomials
$f_1,\ldots,f_r$ in $\Q_p[x]$ such that\begin{enumerate}\item[i)]
$f(x)=f_1(x)f_2(x)\cdots f_r(x),$\item[ii)]the degree of $f_i$ is
$W_i=x_i-x_{i-1},$\item[iii)]all the roots of $f_i$ in
$\overline{\Q}_p$ have $p$-adic valuation $-m_i.$\end{enumerate}
\end{theorem}
\begin{Proof}
See any of the references given above.
\end{Proof}
\begin{corollary}[Dumas]\label{dumas}
With notation as in Theorem \ref{NP},
suppose $f(x)=g(x)h(x)$ is a factorization of $f(x)$ over $\Q_p$.
Then there exist integers $0\leq k_i\leq d_i$ such that $\deg(g) =
\sum_{i=1}^r k_i w_i$.  For each $i=1,\ldots, r$, $f_i$ possesses a
$\Q_p$-irreducible factor of degree at least $w_i$; in particular, $f$
possesses a $\Q_p$-irreducible factor of degree at least $\max(w_1,
\cdots, w_r)$.
\end{corollary}
\begin{Proof}
By the Main Theorem of Newton Polygons, the segments of $NP_p(g)$ and $NP_p(h)$
together make up exactly the
segments of $NP_p(f)$.  Since the $i$th edge of
$NP_p(f)$ is made up of $d_i$ segments of width $w_i=W_i/d_i$, we have
$\deg(g) = \sum_{i=1}^r k_i w_i$ with integers $k_i$ in the range
$0\leq k_i \leq d_i$.  Appealing to the Main Theorem again, we see
that a pure edge must correspond to a $\Q_p$-irreducible
polynomial, giving us the remaining claim.
\end{Proof}
\begin{corollary}[Coleman]\label{corNP}
Suppose $f\in \Q[x]$ and $p$ is a prime.
If an integer $d$ divides the denominator 
(in lowest terms)
of every slope of $NP_p(f)$, then $d$ divides the degree of
any factor $g\in \Q[x]$ of $f$.
\end{corollary}
\begin{Proof}  We give two proofs.  First, this is clearly a special
case of Dumas' corollary (the hypothesis is precisely that each $w_i$
is divisible by $d$).  Now here is Coleman's proof.  By Theorem
\ref{NP}, if $\alpha\in \overline{\Q}_p$ is a root of an irreducible
factor $g$ of $f$, then $p^{\ord_p(n)}$ divides the ramification index
of $\Q_p(\alpha)/\Q_p$ which in turn divides
$[\Q_p(\alpha):\Q_p]=\deg(g)$.  This second proof is a little more
revealing in that it identifies the mechanism behind the divisibility
of the degree of $g$ to be the existence of an inertia group of order
divisible by $d$.
\end{Proof}
\noindent {\bf Remark.} This corollary has in fact appeared a number
of times in the literature, see Mott \cite{mott} and references therein.

Although we will not need it, we mention in passing that the
generalization by Dumas \cite{dumas} of the celebrated Eisenstein
Irreducibility Criterion is a simple consequence of the above
Corollary.
\begin{corollary}[Eisenstein-Dumas Criterion]
Suppose $f=x^n + a_{n-1}x^{n-1} + \cdots + a_1 x + a_0\in \Q[x]$ is
monic polynomial of degree $n$ over $\Q$, and $p$ is a prime.  Let
$m=\ord_p(a_0)$.  Assume $\gcd(m,n)=1$.  If $\ord_p(a_j)\geq m(1-j/n)$
for $j=0,\ldots, n-1,$
then $f$ is irreducible over $\Q$.
\end{corollary}
\begin{Proof}
The geometric meaning of the last hypothesis is that $NP_p(f)$
is ``pure of slope $-m/n$,'' meaning it has only one edge and its slope
is $-m/n$.  Since, by assumption, $\gcd(m,n)=1$, Coleman's Corollary
implies that $n$ divides the degree of any factor in $\Q[x]$ of $f$.
\end{Proof}

Now we recall Coleman's computation of the Newton Polygon of $E_n(x)$
at an arbitrary prime $p$.
Given an integer $n\geq 1$ and a prime $p$, we will define $s+1$
integers $0 = k_0 < k_1 < \cdots < k_s = n$ (where $s$ is the number
of non-zero $p$-adic digits of $n$) called the pivotal indices
associated to $(n,p)$ as follows.  Let us write $n$ in base $p$
recording only the non-zero digits, namely
$$
n = b_1 p^{e_1} + b_2 p^{e_2} + \cdots + b_s p^{e_s}, \qquad 0 < b_1,\ldots
, b_s < p, \qquad e_1 > e_2 > \cdots > e_s \geq 0.$$
The {\em pivotal indices associated to} $(n,p)$ are the partial sums
\begin{equation}
\label{ki}
 k_i = b_1 p^{e_1} + b_2 p^{e_2} + \cdots + b_i p^{e_i}, \qquad i =
0, \ldots, s.
\end{equation}
Note that $k_0=0$ and $k_s=n$.  This definition is motivated by
Coleman's calculation of $NP_p(E_n)$ (see Lemma \ref{ppcc} below).
We will also see that a fundamental fact about the GLP $\LL$ for
$r\geq 0$ is that its $p$-Newton polygons lies on or above
$NP_p(E_n)$.  To explain this, we introduce some more terminology.

\begin{definition}
Suppose $f(x) = \sum_{j=0}^n a_j \frac{x^j}{j!} \in \Q[x]$ and $p$ is
a prime number.  Following P\'olya and Szeg\H o, we call $f$ $p$-{\em
Hurwitz integral} if $\ord_p(a_j)\geq 0$ for $j=0,\ldots,n$.  We call
it {\em Hurwitz integral} if it is $p$-Hurwitz integral for all primes
$p$, i.e. if the Hurwitz coefficients $a_j$ are integral.  We say that
$f$ is  $p$-{\em Coleman integral} if $f$ is
$p$-Hurwitz integral and additionally $\ord_p(a_{k_i}) = 0$ for
$i=0,\ldots,s$ with $k_i$ as defined in (\ref{ki}), i.e. the Hurwitz
coefficients are all $p$-integral and the pivotal
ones are $p$-{\em units}.  
\end{definition}
This definition is motivated by the following Lemma.
\begin{lemma}\label{ppcc}
 If $f\in\Q[x]$ is $p$-Coleman integral of degree $n$, then 
\begin{enumerate}\item[i)]
$NP_p(f)=NP_p(E_n)$;
\item[
ii)] the breaks of $NP_p(f)$ are precisely the pivotal
indices associated to $(n,p)$;
\item[iii)] the slopes of $NP_p(f)$ all have denominator divisible by
$p^{\ord_p(n)}$.
\end{enumerate}
\end{lemma}
\begin{Proof}
We know from Coleman \cite{coleman} that the breaks
of $NP_p(E_n)$ are the pivotal points associated to $(n,p)$.  Since
$f$ is $p$-Hurwitz integral, $NP_p(f)$ lies on or above
$NP_p(E_n)$.  On the other hand, by definition, the corners of
$NP_p(E_n)$ lie on $NP_p(f)$, so $NP_p(f)=NP_p(E_n)$.  The last assertion
iii) follows from ii) and (\ref{ki}).
\end{Proof}

Our proof of Theorem \ref{mt} 
rests on the following two irreducibility criteria.

\begin{lemma}[The Coleman Criterion]
\label{pcc}
Suppose $f\in\Q[x]$ has degree $n$ and $p$ is a prime number.  If $f$
is $p$-Coleman integral, then $p^{\ord_p(n)}$ divides the degree of
any factor $g \in \Q[x]$ of $f$. If $f$ is $p$-Coleman integral for
all primes $p$ dividing $n$, then $f$ is irreducible in $\Q[x]$.
\end{lemma}
\begin{Proof} This is essentially Theorem 1.7 of Sell \cite{wrobel}. 
By Lemma \ref{ppcc}, the
slopes of $NP_p(f)$ all have denominator divisible by $p^{\ord_p(n)}$.
Now apply Corollary \ref{corNP}.
\end{Proof}

Dumas's observed that the Newton Polgyon of the
product of two polynomials is formed by the concatenation, in
ascending slope, of their edges (i.e. is their Minkowski sum, see the
proof of Corollary \ref{dumas}); this is the key tool in the proof of the
following criterion due to Filaseta (see \cite{fschur} for the proof
of a slightly more general version, but note that the convention for
Newton Polygons in that paper differs slighlty from ours).
\begin{lemma}[Filaseta Criterion]
\label{fc}
Suppose $$f(x)=\sum_{j=0}^n b_j \frac{x^j}{j!} \in \Q[x]$$ is
Hurwitz-integral and $|b_0|=1$.  Let $k$ be a positive integer $\leq
n/2$.  Suppose there exists a prime $p\geq k+1$ such that
$$ n(n-1)\cdots (n-k+1) \equiv 0 \bmod{p}, \qquad b_n \not \equiv 0 \bmod{p}.$$
Then $f(x)$ cannot have a factor of degree $k$ in $\Q[x]$.
\end{lemma}

We now give the key calculation allowing the application of the Coleman
Criterion to our family of polynomials.
\begin{lemma}
\label{kl}
i) If $p$ is a prime divisor of $n$, then $\LL$ is $p$-Coleman integral
if and only if $\binom{n+r}{r}\not \equiv 0 \bmod{p}$.

ii) If $\ord_p(n) > \ord_p(r!)$, then $\LL$ is $p$-Coleman integral.
\end{lemma}
\begin{Proof}
From (\ref{lnr}), we see that
$$\LL = \sum_{j=0}^n a_j \frac{x^j}{j!}, \qquad a_j = \frac{
(n-j+1)(n-j+2)\cdots(n-j+r)}{r!},$$ is clearly Hurwitz integral.
From (\ref{ki}) we have $k_0=0$ and we also recall that $$a_0 =
\binom{n+r}{r}=(n+1)\cdots(n+r)/r!.$$  Since $k_i \equiv 0
\bmod{p^{\ord_p(n)}}$ for each $i$, we have $a_{k_i} \equiv a_0
\bmod{p}$.  Thus, the pivotal coefficients $a_{k_i}$ are all $p$-units
if and only if $a_0$ is a $p$-unit, i.e. $\LL$ is $p$-Coleman integral
if and only if $\ord_p(a_0)=0$, proving i)

From the definition of $a_0$, we have 
$$
a_0 \equiv 1 \bmod{p^{\ord_p(n) - \ord_p(r!)}},
$$  
so ii) follows from i).
\end{Proof}
\begin{theorem}
\label{kt}
\begin{enumerate} \item[i)] If
$\gcd(n,\binom{n+r}{r})=1$, then $\LL$ is irreducible over $\Q$.
\item[ii)] If
$\gcd(n,r!)=1$, then $\LL$ is irreducible over $\Q$
\end{enumerate}
\end{theorem}
\begin{Proof}
If $n$ is coprime to $\binom{n+r}{r}$, $\LL$ is
$p$-Coleman integral for every prime divisor $p$ of $n$ by 
Lemma \ref{kl}, so it is irreducible over $\Q$ by the Coleman
Criterion \ref{pcc}.  Part ii), which was first obtained by Sell
\cite{wrobel}, follows from i) since $\gcd(n,r!)=1$ implies
$\gcd(n,\binom{n+r}{r})=1$
\end{Proof}
{\bf Remark.}
In connection with part i) of Lemma \ref{kl}, note that
$p \not |\binom{n+r}{r}$ if and only if there are no ``carries'' in
the addition $n+r$ in base $p$.  Indeed, recalling that
$\ord_p(n!)=\frac{n-\sigma_p(n)}{p-1}$ where $\sigma_p(n)$ is the sum
of the $p$-adic digits of $n$, we have
\begin{eqnarray*}
\ord_p(a_0) &=& \ord_p( (n+r)!) - \ord_p(n!)
- \ord_p(r!)\\
&=& \frac{\sigma_p(n)+\sigma_p(r)-\sigma_p(n+r)}{p-1}.
\end{eqnarray*}
But the latter expression is precisely the number of carries in the base
$p$ addition of $n$ and $r$.  For example, if, say, $n=p$ is prime, then
$\LL$ is irreducible over $\Q$ as long as $-r \not \equiv 1, 2, \cdots, p
\pmod{p^2}$.
More generally, we have
\begin{corollary}\label{cokl}
For each $n$, there is a set of integers $r\geq 0$ of density at least
$\prod_{p|n}
p^{-\ord_p(n)-1}$
for which $\LL$ is irreducible over $\Q$.
\end{corollary}
\begin{Proof}
If $r \equiv 0 \pmod{p^{\ord_p(n)+1}}$, then the addition of $n$ and $r$
in base $p$ cannot have a carry.  Thus, if $r$ is divisible by $\prod_{p|n}
p^{\ord_p(n)+1}$, then by Theorem \ref{kt} and the Remark following it,
$\LL$ is irreducible over $\Q$.
\end{Proof}

\section{Primes in short intervals}

For the proof of Theorem \ref{mt}, we will need to establish the existence
of primes of appropriate size, namely primes for which the Newton
polygon of $\LL$ precludes the existence of factors of certain
degrees.  We will state two such results here, to be used in the next
section.

The first is a well-known consequence of the Prime Number Theorem,
generalizing Chebyshev's Postulate.  For lack of a suitable reference
with an explicit constant, a proof is supplied.

\begin{theorem}\label{ant}
Given $h\geq 2$, there exists a constant $C(h)$ such that 
whenever $N>C(h)$, the interval $[N(1-1/h),N]$ contains a prime.
We may take $$C(h)= e^{h+1/2}(1-1/h)^{-h}.$$
\end{theorem}
\begin{Proof}
We have from Rosser and Schoenfeld \cite{rs}, that
\begin{eqnarray*}
\pi(x) &>& \frac{x}{\log x - 0.5} \qquad \text{for } 67 \leq x \\
\pi(x) &<& \frac{x}{\log x - 1.5} \qquad \text{for } e^{1.5} < x.
\end{eqnarray*}
Since $h\geq 2$, the first inequality applies for $x=N$ and the second
one applies for $x=N-N/h$, assuming only $N\geq 67$.  We then have
\begin{eqnarray*}
\pi(N)- \pi(N-N/h) &>& \frac{N}{\log N -0.5} - \frac{N-N/h}{\log N +
\log(1-1/h) - 1.5}.
%
\end{eqnarray*}
Combining the fractions, the right hand side is positive if and only if
$$
\log N > 1/2 + h - h \log(1-1/h),
$$
proving the lemma, for $N\geq 67$.  We have $C(2)=4e^{2.5}>48$.  
For $N\in [48,67]$,
one easily checks by hand that the lemma holds.
Note that $C(h) \rightarrow e^{h-1/2}$ as
$h\rightarrow \infty$.
\end{Proof}

For Galois group computations in Section 5, we record
\begin{corollary}\label{n+r}
If $n+r\geq 48$ and $n\geq 8+5r/3$, then
there exists a prime $p$ in the interval $(n+r)/2 < p < n-2$.
%
\end{corollary}
\begin{Proof} Apply the Theorem with $h=5$.
\end{Proof}

For the proof of Theorem \ref{mtg}, we will use the following result from
Harborth-Kemnitz \cite{hk},
which is a combination of Theorem \ref{ant} together with a finite
but long computation.
\begin{theorem}[Harborth-Kemnitz]\label{hk}
If $n \geq 48683$, then the interval $(n, 1.001n]$ contains a prime.
\end{theorem}

While Theorem \ref{ant} suffices for the proof of Theorem \ref{mt}, we
may also apply the following stronger, but less concrete, estimate.

\begin{theorem}[Baker-Harman-Pintz \cite{bhp}]\label{bhp}
There is an absolute constant $A$, such that for every $x>A$, the interval
$[x-x^{0.525},x]$ contains a prime.
\end{theorem}

\section{Irreducibility of $\LL$ for large $n$}

We fix $r\geq 0$, and write $n=n_0 n_1=n_2n_3$ where
\begin{equation}
\label{n0}
n_1 =\prod_{p|\gcd(n,\binom{n+r}{r})} p^{\ord_p(n)},
\qquad
n_3=\prod_{\stackrel{p|n}{\ord_p(n) \leq \ord_p(r!)}} p^{\ord_p(n)}.
\end{equation}
Note that the $n_0$ is the largest divisor of $n$ which is coprime to
$\binom{n+r}{r}$.  We also have $n_2 | n_0$ (see the proof of Lemma
\ref{kl}), so $n_1|n_3| \gcd(n,r!)$.  Consequently,
\begin{equation}\label{n1} n_1 \leq
r!,\end{equation}
which is a somewhat crude estimate (see the proof of Theorem \ref{mtg}) but
suffices for the proof of Theorem \ref{mt}.

\begin{lemma} \label{bl}
If there is a prime $p$ 
 satisfying 
$$
 \max(\frac{n+r}{2}, n-n_0)< p \leq n, 
$$
then $\LL$ is irreducible over $\Q$.
\end{lemma}
\begin{Proof}
By Lemma \ref{kl} and Lemma \ref{pcc}, every $\Q[x]$-factor of $f$ has
degree divisible by $n_0$.  If $n_1=1$, then $n=n_0$ and we are done,
so we assume $n_1 > 1$ and proceed by contradiction.  We suppose $f$
has a $\Q[x]$-factor of positive degree $k\leq n/2$.  We know that $$k
\in \{n_0, 2n_0, 3n_0, \ldots, (n_1 -1)n_0\}.$$ To eliminate these
possibilities, we apply the Filaseta Criterion.  Since the latter
requires the constant coefficient to be $1$, we renormalize our
polynomial by setting
\begin{eqnarray*}f(x) &=& a_0^{-1}L^{\langle r \rangle}_n(a_0x)\\
&=& \sum_{j=0}^n b_j \frac{x^j}{j!}
\end{eqnarray*}
with integral Hurwitz coefficients $b_j = a_0^{j-1} a_j$ where $a_0
=\binom{n+r}{r}$.  Note that $b_0=1$ and $b_n = a_0^{n-1}$.  Of
course, the factorization over $\Q$ of $f(x)$ mirrors exactly that of
$\LL$.  With the hypotheses on $p$, we have $p\geq k+1$ (since $k\leq
n/2$).  Moreover, $p\geq n-k+1$ since $k\geq n_0$.  Finally, $p \not |
b_n=a_0^{n-1}$ since $(n+r)/2 <p < n+1$.  Applying the Filaseta
Criterion \ref{fc} to $f(x)$, we find it does not have a factor of
degree $k$, hence neither does $\LL$, giving the desired
contradiction.
\end{Proof}

\begin{lemma}\label{br}
Given $r\geq 0$, there exists a constant $B(r)$ such that 
for every integer $n\geq B(r)$, there exists a prime $p$ satisfying
$$
 \max(\frac{n+r}{2}, n-n_0)< p \leq n,
$$
where $n_0$ is the largest divisor of $n$ coprime to $\binom{n+r}{r}$.
We may take either $$
B(r) = e^{r!+1/2} (1-1/r!)^{-r!} \text{ or } B(r) = \max(A,(r!)^{2.11}),
$$
where $A$ is as in Theorem \ref{bhp}.
\end{lemma}
\begin{Proof}
By (\ref{n1}), $n_1 \leq r!$, so $n-n_0 = n-n/n_1 \leq n(1-1/h)$ with
$h=r!$.  By Lemma \ref{ant}, there exists a prime in the interval
$[n-n_0,n]$ assuming only $n \geq e^{h+1/2}(1-1/h)^{-h}$.  Under this
hypothesis, one easily verifies that $(n+r)/2 < n - n_0$; indeed
merely $n/r>r!/(r!-2)$ suffices.  This establishes the lemma with
$B(r)= e^{r!+1/2}/(1-1/r!)^{r!}$.

Alternatively, if we apply Theorem \ref{bhp} instead, we have $[n -
n_0, n]$ contains a prime once $n > A$ and $n-n/h \leq n -n^{0.525}$,
i.e. if $n> \max(A,(r!)^{2.11})$.  While this gives a better bound
than the one in the previous paragraph (polynomial vs. exponential 
in $r!$), it would be effective only once the constant $A$ is actually
computed.
\end{Proof}

Combining the above Lemmata gives the proof of the first part of
Theorem \ref{mt}.  More precisely, we have proved
\begin{theorem}\label{it}
If $n\geq B(r)$, with $B(r)$ as given in Lemma \ref{br}, then $\LL$
is irreducible over $\Q$.
\end{theorem}

\section{Galois groups}

We begin by recalling a simple criterion based on ramification (as
measured by the Newton polygon) for an irreducible polynomial to have
``large'' Galois group.

\begin{definition}
Given $f\in \Q[x]$, let ${\NN}_f$, called the {\em Newton
Index} of $f$, be the least common multiple of the denominators (in
lowest terms) of all slopes of $NP_p(f)$ as $p$ ranges over all
primes.
\end{definition}
To see that ${\NN}_f$ is well-defined, first note that 0 is
defined to have denominator 1, so slope $0$ segments of $NP_p(f)$
do not contribute to ${\NN}_f$.  On the other hand, for $p$
large enough, all coefficients of $f$ have $p$-adic valuation $0$ so
$NP_p(f)$ consists of a single slope 0 segment.  For a monic
polynomial $f\in \Z[x]$, for example, the Newton Index requires merely the
computation of $NP_p(f)$ for the prime divisors $p$ of its constant
coefficient.  Note also that ${\NN}_f$ divides the least common
multiple of the first $n$ positive integers, where $n=\deg(f)$.

The following result (see Hajir \cite{hajir} for a proof) can be quite
useful for calculating the Galois group of polynomials with
``generic'' ramification.
\begin{theorem}\label{coga}
Given an irreducible polynomial $f \in \Q[x]$, ${\NN}_f$ divides
the order of the Galois group of $f$.  Moreover, if ${\NN}_f$ has a
prime divisor $q$ in the range $n/2 < q < n-2$, where $n$ is the
degree of $f$, then the Galois group of $f$ contains $A_n$.
\end{theorem}

\noindent{\bf Example.}
If $f(x)=L_5^{\langle 3\rangle}(x)$, then $f$ is irreducible over $\Q$
by Lemma \ref{kl}.  An easy calculation shows ${\NN}_f=60$;
indeed we need only consider $p=2,3,5,7$, for which $NP_p(f)$ has
slopes whose denominators are divisible by, repectively, $4,3,5$ and
$2$.  Thus, the Galois group of $f$ has order divisible by $60$.
Since the discriminant of $f$ is not a square (by (\ref{disc}) or
(\ref{disc2}) below), the Galois group of $f$ is $S_5$.
\par

\begin{lemma}
\label{mp}
Suppose $p$ is a prime in the interval $(n+r)/2< p \leq n$.  Then
the $p$-Newton polygon of $\LL$ has $-1/p$ as a slope.
In particular, $p | {\NN}_{\LL}$.
\end{lemma}
\begin{Proof}
Under the assumptions, it is an exercise to calculate the $p$-Newton
polygon of $\LL$ directly from (\ref{lnr}); instead, we use the tools
we have developed to get the result.  According to Lemma \ref{ppcc},
the corners of $NP_p(E_n)$ have $x$-coordinate $0,p$, and $n$ (simply
$0$ and $n$ if $p=n$ of course), so it has $-1/p$ as a slope.  Writing
$\LL= \sum_{j=0}^n a_j x^j/j!$, one checks easily that
$\ord_p(a_0)=\ord_p(a_p)=0$, and we always have $\ord_p(a_n)=0$ since
$a_n=1$.  Since $NP_p(L_n^{\langle r\rangle})$ lies on or above
$NP_p(E_n)$, and they agree at the corners of the latter, they must
coincide.
\end{Proof}

\begin{theorem}\label{gt}
i) If there exists a prime $p$ satisfying $(n+r)/2 < p < n-2$, 
and if $\LL$ is irreducible over $\Q$, then its Galois group over
$\Q$ contains $A_n$.

ii) If $n\geq \max(48-r,8 + 5r/3)$, and if $\LL$ is irreducible over
$\Q$, then its Galois group over $\Q$ contains $A_n$.

iii) For $n > B(r)$ with $B(r)$ as in Lemma \ref{br}, the Galois group
of $\LL$ over $\Q$ contains $A_n$.
\end{theorem}
\begin{Proof} 
We apply Corollary \ref{n+r} in combination with Theorem \ref{coga}
and Lemma \ref{mp}.  For iii), we require Theorem \ref{it} as well.
\end{Proof}
We have thus completed the proof of Theorem \ref{mt}.
We remark that Schur's original method (\cite{sch70}, Satz A), which
was used in \cite{jnt} for the case $r=1$, would yield a proof of
Theorem \ref{gt} as well.
\par
\vspace{\baselineskip}

\noindent {\bf Remark.} By plugging in $\alpha = -1 - n -r$ in Schur's
formula (\ref{disc}), the discriminant of $n!\LL$ is seen to be
\begin{equation}\label{disc2}
\Delta_{n}^{\langle r \rangle} = (-1)^{n(n-1)/2} \prod_{j=1}^{n-1}
(j+1)^{j+1} (r+j)^{n-j}.
\end{equation}
In particular, $\Delta_n^{\langle r \rangle} < 0$, for $n \equiv 2, 3
\pmod{4}$ (recall our blanket assumption $r\geq 0$).  For these values
of $n$, therefore, we know that the Galois group of $\LL$ is not
contained in $A_n$.  If we fix $n>5, n\equiv 0,1 \pmod{4}$, then by
(\ref{disc2}), the Galois group of $\LL$ is contained in $A_n$ if and
only if $r$ is the $x$-coordinate of an integral point on a (fixed)
smooth curve of genus at least $1$, of which there are only finitely
many by Siegel's theorem.  Thus, Conjecture
\ref{galconj} would imply that, for fixed $n$, the Galois group of
$\LL$ is $S_n$ except for a (small) finite number of integers $r\geq
0$.

Similarly, for fixed $r$, if $r$ is small, the proportion of $n$ for
which $\Delta_n^{\langle r \rangle}$ is a square can be large if $r$
is small (as we have already seen for $r=0,1,2$).  
Filaseta has pointed out that this is not so for large $r$.
Specifically, one can check that for $r=3$, $\Delta_n^{\langle r
\rangle}$ is a square if and only if $n\equiv 1 \pmod{4}$ and $n+2$ is
3 times a square; for $r=4,5$, the $n$ for which $\Delta_n^{\langle r
\rangle}$ is a square occur in Fibonacci-type recurrences, namely, for
$r=4$, $n\equiv 0 \bmod{4}$ and $2n+4=\epsilon_3^j+\epsilon_3^{-j}$
for some $j$, and similarly for $r=5$, $n\equiv 1 \bmod{4}$ and $2n+6=
\epsilon_{15}^j+\epsilon_{15}^{-j}$ for some $j$.  Here
$\epsilon_3=2+\sqrt{3}$, $\epsilon_{15}=4+\sqrt{15}$ are the
fundamental units of $\Q(\sqrt{3}),\Q(\sqrt{15})$ respectively.  For
fixed $r\geq 6$, if $n\equiv (r+1)^2 \bmod{4}$, then for $n$ large
enough,$ \Delta_n^{\langle r \rangle}$ cannot be a square because its
$p$-valution must be 1 for some prime $p \in ( (n+r)/2,n+r)$; on the
other hand, if $n\equiv r^2 \bmod{4}$, then integers $n$ for which
$\Delta_n^{\langle r \rangle}$ is a square correspond to integral
points on a smooth curve $y^2 = c_r (x+2) \cdots (x+2\lfloor
r/2\rfloor)$ of positive genus (for some easily determined non-zero
constant $c_r$); there are, therefore, only finitely
many such $n$ by Siegel's theorem.

\section{Properties of $\LL$ for $n\leq 4$}

In this section, as well as the next, we establish more evidence for
Conjectures \ref{conj} and \ref{galconj} of a somewhat complementary
nature to Theorem \ref{mt}.  Namely, we fix $n$ and consider those
$\alpha\in \Q$ for which $L_n^{(\alpha)}(x)$ is irreducible over $\Q$.
This point of view has a rather different flavor.  For arbitray $n$, 
the methods of this
paper allowed us to get only a weak result (Corollary \ref{cokl}) in
this direction.  If $n\geq 5$, a much more fruitful,
algebro-geometric, point of view, adopted in \cite{HW}, is to consider
the covering of curves $\X_1 \rightarrow \P^1$ given by the
projection-to-$y$ map, where $\X_1:{\mathcal L}_n^{(y)}(x)=0$ is the
projective curve defined by the $n$th degree GLP. The Galois closure
of this cover, call it $\X'$, has monodromy group $S_n$ (by Schur's
result that ${\mathcal L}_n^{(0)}(x)$ has Galois group $S_n$).  By
estimating from below the genus of $\X_1$ and other quotients of
$\X'$, the following theorem was proved in \cite{HW}.

\begin{theorem}[Hajir-Wong]\label{hajir-wong}
Suppose an integer $n\geq 5$ and a number field $K$ are fixed.  There
is a finite subset ${\mathcal E}(n,K)\subset K$ such that for $\alpha
\in K - {\mathcal E}(n,K)$, we have i) $L_n^{(\alpha)}(x)$ is
irreducible over $K$, and ii) the Galois group of $L_n^{(\alpha)}(x)$
contains $A_n$ (if $5\leq n \leq 9$), is the full symmetric group (if
$n\geq 10$).  
\end{theorem}

Applying the theorem with $K=\Q$, we have the following
nice complement to the main theorem \ref{mt} of this paper.

\begin{corollary}
For each $n\geq 5$, there is a bound $C_n$ such
that 
Conjectures \ref{conj} and \ref{galconj} hold for the pair $(n,r)$ whenever
$r\geq C_n$.
\end{corollary}
\noindent {\bf Remark.}
The constant $C_n$ in the above Corollary is ineffective since the
proof of the Theorem preceding it rests on Faltings' theorem on
finitude of rational points on curves of genus at least 2; for the
Corollary, we could apply Siegel's theorem on integral points instead,
but this does not resolve the effectivity issue either since for
$n\geq 5$, the relevant curves have genus greater than 1.

For $n\leq 4$, on the other hand, GLP admitting proper factors over
$\Q$ turn out to be plentiful, as such factors correspond to rational
points on certain curves of genus 0 or 1.  In this section, we
calculate the (very few) {\em integral} points on these curves
effectively, thereby establishing Conjectures \ref{conj} and
\ref{galconj} for $n\leq 4$ and all $r\geq 0$.  We summarize the
results in the following theorem.  During the proof, we will give
parametrizations for all $\alpha\in \Q$, $n\leq 4$, for which
$L_n^{(\alpha)}(x)$ is $\Q$-reducible.  We also parametrize, for
$n=4$, an infinite family of specializations which are reducible but
have exceptional Galois group $D_4$.

\begin{theorem}\label{n3}
(a) If $n\leq 4$ and $r\geq 0$, then $\LL$ is irreducible over $\Q$ and
has Galois group containing $A_n$.  If $n\leq 3$, this Galois group is
in fact the full symmetric group $S_n$.  

(b) For each $n \in \{ 2,3,4\}$, there exist infinitely many rational
numbers $\alpha$ such that $L_n^{(\alpha)}(x)$ is reducible over $\Q$.

(c) There are infinitely many rational numbers $\alpha$ for which
$L_4^{(\alpha)}(x)$ is irreducible over $\Q$ with Galois group {\em
not} containing $A_4$.
\end{theorem}
\begin{Proof}
To prove irreducibility of $\LL$ for a fixed $n$, and arbitrary
$r\geq 0$, the techniques we have used so far (the existence of
ramification at primes dividing $n!$) would have to be modifed,
because for suitable $r$, not all primes less than $n$ ramify in
the splitting field of $L_n^{\langle r \rangle}(x)$ over $\Q$.
We can take a more direct approach. 
For $n=2$, the sign in the discriminant formula (\ref{disc2}) is already enough
to show the irreducibility of all $\LL$ for $n=2$, and the same
formula shows that $L_2^{(\alpha)}(x)$ is reducible exactly when
$\alpha+2$ is a rational square.  It also shows that
$L_3^{\langle r \rangle}(x)$ does not have Galois group $A_3$.
 
Now suppose $n=3$.
Let $s=r+1$ and put 
$$
f(x) := 3! L_3^{\langle r \rangle}(x-r-1) = x^3 + 3 s x + 2 s.
$$ We need to show that $f(x)$ is irreducible over $\Q$.  It suffices
to show that $f$ does not vanish on $\Z$.  Suppose $f(m)=0$ for some
integer $m$.  Writing $$s=\frac{-m^3}{3m+2},$$ we see that for an odd
prime $p$ dividing $s$, $\ord_p(s)=3\ord_p(m)$ because $p|s$ implies
$p|m$ which implies $p\not| 3m+2$.  Let us write $s=2^as_0$,
$m=2^bm_0$ where $s_0$ and $m_0$ are odd integers.  We then have
\[\label{3b} 2^{b-1}\cdot 3\cdot m_0 + 1 = -2^{3b-a-1}.\] Thus,
$3b\geq a+1$.  If $3b=a+1$, then $2^{b-1}\cdot 3\cdot m_0=-2$ which is
not possible, so $3b>a+1$.  By (\ref{3b}), $2^{b-1}\cdot 3 \cdot m+1$
is even, so we must have $b=1$.  But then $a \in \{0,1,2\}$ and each
of these is easily eliminated.  Thus, $f(x)$ is irreducible over $\Q$.
Moreover, we see immediately that
$L_3^{(\alpha)}(x)$ is {\em reducible} over $\Q$ for infinitely many
rational numbers $\alpha$, and that this is so exactly for those of the form
$$\alpha=\frac{m^3-9m-6}{3m+2}, \qquad m\in \Z.$$

For $n=4$, we consider linear factors and quadratic factors separately.
We start by simplifying the model via killing the trace term as before,
i.e. we reparametrize with $s=r+1$ again and define
$$ g(x,s) := 4! L_4^{\langle s-1 \rangle}(x-s) = x^4 + 6 s x^2 + 8 s x
+3 s^2+6s.
$$ A $\Q$-linear factor $(x-x_0)$ of $g(x,s_0)$ for a rational number
$s_0$ corresponds exactly to a (finite) rational point $(x_0,s_0)$ on
the curve $\X_1: g(x,s)=0$.  It is easy to see that this curve has
genus 1, so is elliptic ($(0,0)$ is on it).  Upon using the
Cayley-Hermite formula, (implemented in Maple 7 for
example), to put ${\mathcal X}_1$ in Weierstrass form, we find it is
birational to the minimal Weierstrass model $384H2: Y^2 = X^3 + X^2 -
25X + 119$, of conductor $384=2^7\cdot 3$, where 
$$ x =6\frac{4X+Y+28}{X^2-22X-95},
 \qquad s=-216\frac{X^2+10X+8Y+129}{X^4-44X^3+294X^2+4180X+9025}
.$$
Here we are using the notation from Cremona's table (available, for
instance, in a very usable format at \cite{cremona}), from which we
learn that this elliptic curve has infinite Mordell-Weil group over
$\Q$, generated by the point $P_1=(-1,12)$ of infinite order and the
$2$-torsion point $P_0=(-7,0)$.  This completes the proof of (b).
By the usual height arguments, it is
not difficult to show that the only integral points on $g(x,s)=0$ are
$$(0,0),(0,-2),(3,-1),(4,-2),(-1,-1),(-2,-2),(-3,-3), (3,-27),
(-3,-9).$$
All but the last two of these 
correspond to the trivial
factorizations (see (\ref{factor})).  This verifies that for $n=4$ and
integers  $r\geq 0$ (as well as integers $r\leq -11$), 
$\LL$ does not have a linear
factor over $\Q$.  Note the exceptional factorization for $s=-9,-27$, i.e.
$r=-10,-28$,
corresponds to the factors $x-6$ and $x-30$ in $L_4^{(5)}(x)$ and 
$L_4^{(23)}(x)$ respectively.

The quadratic factors of $L_4^{\langle r \rangle}(x)$ are also
parametrized by a curve ($\X_2$ let us call it), for which we can find
a model by writing $$ g(x,s) = x^4+6sx^2+8sx+3s^2+6s= (x^2 +Ax
+B)(x^2-Ax+C)$$ and equating coefficients.  A simple elimination of
the resulting equations gives us the curve
$$\X_2: h(A^2,s)=0,$$
where
$$h(z,s):= z^3 + 12sz^2+24s(s-1)z-64s^2,$$
is the cubic resolvent of $g(x,s)$.  One checks that $\X_2$
also has genus $1$ and is birational to $384H1:Y^2 = X^3+X^2-35X+69$
via 
$$A = \frac{-6Y}{X^2+4X-23}, 
\qquad 
s = \frac{-27(X-3)^2}{(2X-5)(X^2+4X-23)}.
$$ Thus, $\X_1$ and $\X_2$ in fact form an isgoney class of order 2.
(Note in passing that, with respect to the projection-to-$s$ map, 
the fiber product $\X'=\X_1 \times_{{\P^1}}
\X_2$ is the minimal Galois cover of either).  In particular, $X_2$
also has rank $1$, with Mordell-Weil group generated by $P_1=(1,6)$
together with $2$-torsion point $(3,0)$.  We find the integral points
on this curve correspond exactly to the previously known trivial
factorizations, namely $(0,0),(\pm 2,-2),(\pm
2,-1),(\pm 4,-2)$.  This completes the proof of conjecture \ref{conj}
for $n\leq 4$.  

Turning to the Galois group over $\Q$ of $g(x,s)$, we
know that it contains $A_4$ if and only if the cubic resolvent
$h(z,s)$ does not have a rational root, i.e. if and only if the curve
$\Y_2:h(z,s)=0$, over which $\X_2$ is a double cover, does not have
a $\Q$-rational point. 
Considering $h(z,s)$ as a quadratic in $s$
with discriminant $(4z)^2(3z^2-20z+36)$, we see that the integral (or
rational) points on $\Y_2$ correspond the integral (rational) points
on the conic $w^2 = 3z^2-20z+36$.  This already suffices to prove (c),
and one can give an explicit formula
$$
s= \frac{z(12-6z\pm \sqrt{3z^2-20z+36})}{8(3z-8)}, \qquad (3z-10)^2-3w^2=-8,
$$
for rational values of the parameter $s$ at which $L_4^{\langle
s-1\rangle}(x)$ has dihedral Galois group $D_4$ (hence is not
contained in $A_4$); it is clear that the values of $s,w,z$ can be
parametrized by the trace of powers of the fundamental unit of
$Z[\sqrt{3}]$ or a corresponding suitable recurrence.  If $s$ is
restricted to the integers, then by Gauss's Lemma, $z$ and $w$ are
also integers, and one again shows that $s=0,-1,-2$ give the only
integral points on the model $\Y_2$; we omit the details.  This
completes the proof of Conjecture \ref{galconj} for $n\leq 4$, as well
as that of the theorem.
\end{Proof}

\noindent{\bf Remark.}  The Galois group of $L_4^{\langle 4
\rangle}(x)$ is $A_4$ for infinitely many integers $r$, namely exactly
those expressible as $r=-2+\sqrt{12k^2+1}$, with $k\in \Z$ (these can
be parametrized by the trace of powers of the fundamental unit of
$\Z[\sqrt{3}]$), or by a suitable recurrence.

\section{Proof of Theorem \ref{mtg}}

Now we want to prove Conjectures \ref{conj} and \ref{galconj}
for arbitrary $n$ and small $r$.
\begin{proof}[Proof of Theorem \ref{mtg}.] 
As mentioned earlier, the cases $r=0,1,2$ have already 
appeared in the literature, missing only the calculation
of a few Galois groups for small $n$.  Since it is no extra work
we give a uniform proof here for all $0\leq r \leq 8$.  By
Theorem \ref{it}, this has been reduced to a finite
calculation, but the bound given there is prohibitively
large, since $B(8)$ is greater than $2\cdot 10^{17511}$.

We begin by sharpening the bound (\ref{n1}).
Recall our notation that $n_0$ is the largest divisor of $n$ 
coprime to $\binom{n+r}{n}$, and $n_1=n/n_0$ is its 
complement.  We have $n_1 | \gcd(n,r!)$.
 
We claim that for $r\leq 8$ and all $n\geq 1$, $n_1 \leq 840$.  We
know that in this range, $n_1 | 8!=2^7\cdot3^2\cdot 5\cdot 7$, so it
suffices to prove that $\ord_2(n_1) < 4$, $\ord_3(n_1) < 2$.  Both of
these facts follow easily from the following observation.  Recall that
a prime $p$ divides $\binom{n+r}{n}$ if and only if there is a carry
in the base-$p$ addition of $n$ and $r$.  Thus, if $n\equiv 0
\pmod{p^a}$, and $r< p^a$, then $p$ does not divide $\binom{n+r}{n}$
so $p$ does not divide $n_1$.  Since $8 < 2^4$ and $8<3^2$, we are done.
In general, by this argument we have, for a given fixed $r$, that
$$
n_1 \leq \prod_{p|r!} p^{\lfloor \log_p(r)\rfloor}.
$$

Thus, for $0\leq r \leq 8$ and $n\geq 1$, we have
$n-n_0=n(1-1/n_1)\leq (839/840)n$.  By Theorem \ref{hk}, the interval
$(839n/840,n]$ contains a prime for $n\geq 48742$ (note that
$1/839=0.00119\ldots > 0.001$; one checks easily then that we can
replace 48742 by 44350 if we wish).  For $0\leq r \leq 8$, $n\geq 9$,
we have $n-n_0\geq (n+r)/2$; we have therefore shown that for $0\leq r
\leq 8$, $n\geq 48742$, there exists a prime $p$ in the range $
\max( (n+r)/2,n-n_0) < p \leq n$.  This proves the irreducibility of
$\LL$ for $n\geq 48742$ by Lemma \ref{bl}.

Now we need to handle the small degrees.
By Theorem \ref{n3}, we can take $n\geq 4$.  Using PARI, for each pair
$(n,r)$ in the box $4\leq n \leq 48741$, $0\leq r \leq 8$, we
calculated $n_0$ and checked 
i) whether $n=n_0$, and ii) whether the smallest prime exceeding
$\max((n+r)/2,n-n_0)$ is at most $n$
(PARI is equipped with a table of primes).
If i) holds, then $\LL$ is
irreducible by Theorem \ref{kt}, and if ii) holds, then
$\LL$ is irreducible by Lemma \ref{bl}.  It took PARI
only a few seconds to verify that among these 438642 pairs $(r,n)$,
only 24 cases remain (listed in Table 1) where neither Lemma \ref{bl}
nor Theorem \ref{kt} applies.  We verified using PARI's routine
\verb?polisirreducible? that for these remaining pairs, $\LL$ is
irreducible.  

In order to supply a more tangible certificate of irreducibility, we
list in Table 1, with one exception, a prime $\l$ such that the
reduction $\LL$ is irreducible in $\F_\l[x]$.  The pair $(4,5)$ is
exceptional because the discriminant of $L_4^{\langle 5 \rangle}(x)$
is a square, so by a theorem of Stickelberger, this polynomial is never
irreducible over a prime field $\F_\l$.  It is simple enough to check
that $L_4^{\langle 5 \rangle}(x)$ has no linear factor, and we can
verify that it has no quadratic factor, for example, by applying Lemma
1 from \cite{fl} to $4!L_4^{\langle 5\rangle}(x)$ with
$k=2,\ell=1,p=7$.  The very last entry in the table is also
interesting.  Although $L_{120}^{\langle 8 \rangle}$ is not
$p$-Coleman integral for any prime divisor $p$ of $120$, one checks
that all slopes of its $p$-Newton polygon are divisible by $p$ for
$p=3$ and $p=5$.  By Corollary \ref{corNP},
$15$ divides the degree of any irreducible factor of
$L_{120}^{\langle 8 \rangle}$.  Thus, even though $n_0=1$, we can
apply Lemma \ref{kt} with $n_0=15$ and $p=107$ to get the
irreducibility of $ L_{120}^{\langle 8 \rangle}$.

Now let us turn to the computation of the Galois group.  Of course, we
need only consider $n\geq 4$.  When $n<8$, $(n/2,n-2)$ does not
contain prime, so we cannot apply Jordan's criterion.  For the 36
polynomials $\LL$ with $0\leq r \leq 8$ and $4\leq n \leq 7$, we used
the PARI routine \verb?polgalois? to verify that the Galois group
contains $A_n$.  

Now suppose $n\geq 8$ and $r\leq 8$.  By Theorem \ref{gt}, ii), we are
done if $n\geq 48$.  Of the remaining pairs $(r,n)$, when
$((n+r)/2,n-2)$ contains a prime, we apply Theorem \ref{gt}, i).
There remain 47 cases, listed in Table 2.  In these 47 cases, since
$n\geq 8$, there exists a prime in the interval $(n/2,n-2)$, labelled
$q$ in Table 2.  We check in each case that $NP_q(\LL)$ has at least
one slope with denominator $q$, then apply Theorem \ref{coga}.  Thus, in all
cases, the Galois group of $\LL$ contains $A_n$.
\end{proof}


\begin{center}
\begin{tabular}{|c|c|c|}
\hline
$r$ & $n$ & $\l$ \\
\hline
$3 $ & $ 6 $ & $ 13 $\\ \hline
$4 $ & $ 4 $ & $ 17 $\\ \hline
$4 $ & $6  $ & $ 29  $ \\ \hline
$5 $ & $4  $ & $*  $ \\ \hline
$5 $ & $6  $ & $ 23 $  \\ \hline
$5 $ & $20  $ & $  149$ \\ \hline
$6 $ & $4  $ & $ 13 $ \\ \hline
$6 $ & $6  $ & $31  $ \\ \hline
\end{tabular}
\begin{tabular}{|c|c|c|}
\hline
$r$ & $n$ & $\l$ \\
\hline
$6 $ & $10  $ & $ 17 $ \\ \hline
$6 $ & $12  $ & $ 29 $ \\ \hline
$6 $ & $20  $ & $ 311 $ \\ \hline
$7 $ & $4  $ & $ 13 $ \\ \hline
$7 $ & $6  $ & $ 47 $ \\ \hline
$7 $ & $10  $ & $ 47 $ \\ \hline
$7 $ & $12  $ & $ 47 $ \\ \hline
$7 $ & $20  $ & $ 271 $ \\ \hline
\end{tabular}
\begin{tabular}{|c|c|c|}
\hline
$r$ & $n$ & $\l$ \\
\hline
$7 $ & $42  $ & $  79$ \\ \hline
$8 $ & $6  $ & $ 17 $ \\ \hline
$8 $ & $8 $ & $29  $ \\ \hline
$8 $ & $10  $ & $137  $ \\ \hline
$8 $ & $12  $ & $173  $ \\ \hline
$8 $ & $24  $ & $191  $ \\ \hline
$8 $ & $42  $ & $113  $ \\ \hline
$8 $ & $120  $ & $613  $ \\
\hline
\end{tabular}
\end{center}

\begin{center}
{\bf Table 1}
\end{center}

\begin{center}
\begin{tabular}{|c|c|c|}
\hline
$r$ & $n$ & $q$ \\
\hline
$ 1$ & $ 9$ & $ 5$ \\ \hline
$ 1$ & $ 13$ & $ 7$ \\ \hline
$ 2$ & $ 8$ & $ 5$ \\ \hline
$ 2$ & $ 9$ & $ 5$ \\ \hline
$ 2$ & $ 12$ & $ 7$ \\ \hline
$ 2$ & $ 13$ & $ 7$ \\ \hline
$ 3$ & $ 8$ & $ 5$ \\ \hline
$ 3$ & $ 9$ & $ 5$ \\ \hline
$ 3$ & $ 11$ & $ 7$ \\ \hline
$ 3$ & $ 12$ & $ 7$ \\ \hline
$ 3$ & $ 13$ & $ 7$ \\ \hline
$ 4$ & $ 8$ & $ 5$ \\ \hline
$ 4$ & $ 9$ & $ 5$ \\ \hline
$ 4$ & $ 10$ & $ 7$ \\ \hline
$ 4$ & $ 11$ & $ 7$ \\ \hline
$ 4$ & $ 12$ & $ 7$ \\ \hline
\end{tabular}
\begin{tabular}{|c|c|c|}
\hline
$r$ & $n$ & $q$ \\
\hline
$ 4$ & $ 13$ & $ 7$ \\ \hline
$ 5$ & $ 8$ & $ 5$ \\ \hline
$ 5$ & $ 9$ & $ 5$ \\ \hline
$ 5$ & $ 10$ & $ 7$ \\ \hline
$ 5$ & $ 11$ & $ 7$ \\ \hline
$ 5$ & $ 12$ & $ 7$ \\ \hline
$ 5$ & $ 13$ & $ 7$ \\ \hline
$ 6$ & $ 8$ & $ 5$ \\ \hline
$ 6$ & $ 9$ & $ 5$ \\ \hline
$ 6$ & $ 10$ & $ 7$ \\ \hline
$ 6$ & $ 11$ & $ 7$ \\ \hline
$ 6$ & $ 12$ & $ 7$ \\ \hline
$ 6$ & $ 13$ & $ 7$ \\ \hline
$ 7$ & $ 8$ & $ 5$ \\ \hline
$ 7$ & $ 9$ & $ 5$ \\ \hline
$ 7$ & $ 10$ & $ 7$ \\ \hline
\end{tabular}
\begin{tabular}{|c|c|c|}
\hline
$r$ & $n$ & $q$ \\
\hline
$ 7$ & $ 11$ & $ 7$ \\ \hline
$ 7$ & $ 12$ & $ 7$ \\ \hline
$ 7$ & $ 13$ & $ 7$ \\ \hline
$ 7$ & $ 15$ & $ 11$ \\ \hline
$ 7$ & $ 19$ & $ 11$ \\ \hline
$ 8$ & $ 8$ & $ 5$ \\ \hline
$ 8$ & $ 9$ & $ 5$ \\ \hline
$ 8$ & $ 10$ & $ 7$ \\ \hline
$ 8$ & $ 11$ & $ 7$ \\ \hline
$ 8$ & $ 12$ & $ 7$ \\ \hline
$ 8$ & $ 13$ & $ 7$ \\ \hline
$ 8$ & $ 14$ & $ 11$ \\ \hline
$ 8$ & $ 15$ & $ 11$ \\ \hline 
$ 8$ & $ 18$ & $ 11$ \\ \hline 
$ 8$ & $ 19$ & $ 11$ \\ \hline
$~$ & $~$ & $~$ \\
\hline
\end{tabular}
\end{center}
\begin{center}
{\bf Table 2}
\end{center}
\section{A Question}

Given $f(x)=\sum_{j=0}^n a_j x^j \in \Q[x]$, let us say
$g(x)=\sum_{j=0}^n a_j b_j x^j$ is an {\em admissible modification} of
$f(x)$ if $b_j \in \Z$ for all $0\leq j\leq n$ and $b_0=\pm 1, b_n=1$.
We could also allow $b_n=-1$, but since multiplication by $-1$ is
harmless when it comes to irreducibility and Galois groups, we can
dispense with it.

Already in Schur's original treatment of $E_n(x)=\L{0}$, he proved not
just the irreducibility of $E_n$ but also of all its admissible
modifications.  In \cite{ft}, Filaseta and Trifonov prove the
irreducibility of all admissible modifications of the Bessel
polynomials $z_n(x) = n!\L{n}$.  Also, the Filaseta-Lam theorem quoted
in the introduction was in fact proved for all admissible
modifications of $L_n^{(\alpha)}(x)$ for $n$ large enough with respect
to $\alpha$.  These results, combined with Conjecture \ref{conj}
suggest the following question.
\begin{question}
For which pairs of non-negative integers $(r,n)$ is it true that
every admissible modification of $\LL$ is irreducible over $\Q$?
\end{question}

The particular strategy developed in this paper would not appear to be
suitable for answering this question, but techniques of \cite{ft} and
\cite{fl}, suitably altered, would hopefully apply.

Some experimentation reveals that there are exceptions already for
$n=2$.  Indeed, suppose $r=4m^2-1$ and the modifying coefficients
$(b_0, b_1, b_2)$ are $(-1,m,1)$.  The resulting admissible
modification of $2L_2^{\langle r \rangle}(x)$ is
$$
x^2 + 8m^3x-4m^2(4m^2+1) = (x-2m)(x+2m+8m^3).
$$
If one does not allow the modification of the constant coefficient,
then it is not hard to show that the resulting admissible modifications
of $L_2^{\langle r \rangle}(x)$ are always irreducible over $\Q$.
Moreover, a PARI calculation for $n=3$ and $r\leq 100$, 
with modification
coefficients $(b_0, b_1, b_2, b_3)$ satisfying $|b_0|=1$, $b_3=1$,
$|b_1|, |b_2| \leq 100$ turned up only irreducible polynomials (more than
2 million of them).

\vspace{1in}

\begin{sc}
\noindent Farshid Hajir\\
 Dept. of Mathematics \& Statistics\\
 University of Massachusetts, Amherst\\
 Amherst MA 01003 \\
\end{sc}
 \verb?hajir@math.umass.edu?

\end{document}